\documentclass[12 pt]{amsart}
\usepackage{amscd,amssymb}
\usepackage[arrow,matrix]{xy}
\usepackage{graphicx}
\oddsidemargin=0.3in \evensidemargin=0.3in \theoremstyle{plain}
\newtheorem{thm}[subsection]{Theorem}
\newtheorem{lem}[subsection]{Lemma}
\newtheorem{prop}[subsection]{Proposition}
\newtheorem{cor}[subsection]{Corollary}

\theoremstyle{definition}
\newtheorem{rk}[subsection]{Remark}
\newtheorem{defn}[subsection]{Definition}
\newtheorem{ex}[subsection]{Example}

\numberwithin{equation}{section} 

\newcommand{\Z}{\mathbb{Z}}

\newcommand{\C}{\mathbb{C}}

\newcommand{\PP}{\mathbb{P}}

\DeclareMathOperator{\Hom}{Hom} 
 
\DeclareMathOperator{\bideg}{bideg}

\DeclareMathOperator{\multideg}{multideg}

\begin{document}

\title [Topology and Factorization of Polynomials]
{Topology and Factorization of Polynomials}

\author[Hani Shaker]{Hani Shaker }
\address{ School of Mathematical Sciences, GCU, 68-B New Muslim Town, Lahore Pakistan.}
\email {hani.uet@gmail.com}


\begin{abstract}

For any polynomial $P \in \C[X_1,X_2,...,X_n]$, we describe a $\C$-vector space
$F(P)$ of solutions of a linear system of equations coming from some algebraic partial differential equations such that the dimension of $F(P)$
is the number of irreducible factors of $P$. Moreover, the knowledge of $F(P)$ gives a complete factorization
of the polynomial $P$ by taking gcd's.
This generalizes previous results by Ruppert and Gao
in the case $n=2$.

\end{abstract}
\maketitle


\section{Introduction}

Let $K$ be the algebraic closure of a field $k$ and
let $k[X_1,X_2,...,X_n]$ be the polynomial ring in $n$ indeterminates.
The zero set of
a polynomial $P \in k[X_1,X_2,...,X_n]$ of $\deg d >0$
is a hypersurface $V(P)$ in $K^n.$
As the polynomial ring is a factorial ring, we can write
$P=\prod_{i=1}^sP_i,$ where $P_i$ are the irreducible factors of $P$ in
$K[X_1,X_2,...,X_n].$

We assume that the factors $P_i$ are distinct,
i.e. $P$ is a {\it reduced polynomial}. The prime factorization of $P$
corresponds to the decomposition into irreducible components
$V(P)=\prod_{i=1}^sV(P_i)$ of the hypersurface $V(P)$.

\medskip

A natural question to ask is: {\it "How can we compute $s$, the number of irreducible factors  of $P$ (resp. irreducible components of $V(P)$)
from the coefficients of $P$?"}
A variant of this problem (called the {\it absolute factorization problem})
is when $P$ is assumed to be irreducible in $k[X_1,X_2,...,X_n]$, see  \cite{GC}.

\medskip

In this paper we recall in Section 2 briefly Gao's results in the
case $n=2$, see \cite{Gao}, and then some usual techniques for
reducing the case $n>2$ to the case $n=2$ by taking generic linear
sections, see \cite{GC}.

Since all these reduction techniques are not easy to use in practice
(since the notion of a {\it generic} linear section is quite subtle
as we show by some examples), we develop in Sections 3 and 4 of our
note a direct approach to the case $n>2$.

Unlike Ruppert-Gao's approach, which is purely algebraic and works over any
field $k$ of characteristic zero or of relatively large characteristic, our approach is topological, using de Rham cohomology, and hence works
only for the algebraically closed subfields of the field of complex numbers $\C$.

\section{Ruppert-Gao's idea and the reduction techniques}

Assume that $n=2$ and denote by $X,Y$ the two indeterminates.
If $P(X,Y)={\prod_{i=1}^sP_i(X,Y),}$ is the factorization of $P$ into irreducible factors in $K[X,Y]$,
then by taking the
partial derivatives on both sides, we have
\begin{equation} \label{eq1}
P_X=\sum_{i=1}^s(\prod_{j\neq i}P_j)
\frac{\partial P_i}{\partial X}=\sum_{i=1}^sg_i \text{ where } g_i=(\prod_{j\neq i}P_j) \frac{\partial
P_i}{\partial X}
\end{equation}

and also
\begin{equation} \label{eq2}
P_Y=\sum_{i=1}^s(\prod_{j\neq i}P_j)
\frac{\partial P_i}{\partial Y}=\sum_{i=1}^sh_i\text{ where } h_i=(\prod_{j\neq i}P_j) \frac{\partial
P_i}{\partial Y}.
\end{equation}

Note that we can write
$$\frac{\partial}{\partial X}(\log
P_i)=\frac{1}{P_i}\frac{\partial P_i}{\partial X}=\frac{g_i}{P}, \ \ \ \frac{\partial}{\partial Y}(\log
P_i)=\frac{1}{P_i}\frac{\partial P_i}{\partial Y}=\frac{h_i}{P}$$
which yields
\begin{equation} \label{eq3}
\frac{\partial}{\partial Y}(\frac{g_i}{P})=\frac{\partial}{\partial X}
(\frac{h_i}{P})\ \ \ \ \ \ \ \ \ \mbox{for}\ \ \ i=1,...,s. \ \ \ \
\end{equation}

\begin{defn} \label{d1}
Let $P(X,Y)\in K[X,Y]$ be such that ${\deg}_X(P)=m_1$ , $
{\deg}_Y(P)=m_2$. Then the bidegree of $P$ is defined as
$${\bideg}(P)=({\deg}_X(P),{\deg}_Y(P))=(m_1,m_2).$$
\end{defn}
In our case, we obviously have
\begin{equation}
\nonumber {\bideg}(g_i) \leq (m_1-1,m_2) \text{ and } {\bideg}(h_i) \leq (m_1,m_2-1).
\end{equation}
\begin{defn} \label{d2}
Let $F(P)$ be the $K-$vector space of solutions
$(v,w)\in K[X,Y]^2$ of the partial differential equation
$${\frac{\partial}{\partial Y}(\frac{v}{P})=\frac{\partial}{\partial
X} (\frac{w}{P})}$$
 such that ${{\bideg}(v)\leq(m_1-1,m_2)}
,{{\bideg}(w)\leq(m_1,m_2-1).}$
\end{defn}
This partial differential equation was first considered by Ruppert \cite{R1}, \cite{R2}. Moreover, it was clear to Ruppert and Gao that this is just the condition that a certain 1-form is closed, see the comment just before Theorem
2.1 in \cite{Gao}.

\begin{thm}\emph{\textbf{(Gao's Theorem)}}\label{t1}\ \cite{Gao}\\
If $P(X,Y)={\prod_{i=1}^sP_i(X,Y)}$ is the factorization of $P$ into irreducible factors in $K[X,Y]$, then
$s=\dim_KF(P)$  and the set
$${\{(g_i,h_i)\ |\ i=1,...,s\}}$$
is a basis for $F(P)$.
\end{thm}
\begin{cor} \label{c1} \ \\
\noindent (i) $P$ is irreducible if and only if
$\dim_KF(P)=1.$

\noindent (ii) $P_i=\mbox{g.c.d.}(P,v-\lambda_iP_X)$  where
$v=\Sigma_{i=1}^s\lambda_ig_i$ is a generic vector in the vector
space $E(P)$ obtained from $F(P)$ by projecting on the first factor.
\end{cor}
Here $\lambda=(\lambda_1,...,\lambda_s) \in \C^s$ and the genericity means that $\lambda$ has to avoid a proper Zariski closed subset of $\C^s$.
The first claim is an obvious consequence of Theorem \ref{t1} and was obtained already by Ruppert \cite{R1}. The second one is much more subtle
and we will discuss this point in the general case in the last section, see in particular Proposition
\ref{p6}.

\medskip

Now we return to the general case $n \ge 2$.
Let $V(P)$ be the affine hypersurface defined by $P=0$ in the affine space $K^n$.
Let $E$ be an affine
plane in $K^n$ such that $ V(P) \cap E$ is a curve in $E$.
One may ask {\it "Is
there some relation between the number of irreducible components of
$V(P)$ and $V(P)\cap E$?} or, more precisely: {\it Are these numbers always  equal?"}
The answer is to such questions depends
on the choice of $E.$ Let us look at two examples.

\begin{ex} \label{e1}

\noindent (i) Consider the Whitney umbrella $S:x^2-zy^2=0$, an irreducible singular surface in $\C^3$. Choose two
planes
$E_0:z=1$ and $E_1:y=1$. One can see that $S\cap
E_0$ is the union of two lines, namely $x^2-y^2=0$, and  $S\cap E_1$
is irreducible and isomorphic to $\C$.

\noindent (ii) Consider the smooth irreducible surface $S': x^2y-x-z=0$ in $\C^3$. Choose two
planes
$E_0:z=0$ and $E_1:z=1$. One can see that $S'\cap
E_0$ has two components $x=0$ and $xy-1=0$, while  $S' \cap E_1$ is irreducible, and isomorphic to $\C^*$.

\end{ex}

By Bertini's second Theorem we know that the number of irreducible components of
$V(P)$ and of $V(P)\cap E$ coincide if the the 2-plane $E$ is generic, see \cite{GC}, subsection 9.1.3 for
an excellent survey of this problem as well as Section 5 in \cite{Gao}, for relations to an effective Hilbert irreducibility theorem. In practice it is quite difficult to decide when a given plane $E$ is
generic. In the next section we explain the relation between this genericity and transversality to some
Whitney regular stratifications, but this is not easy to check on explicit examples.

Moreover, once we have the factorization of $P$ in the plane $E$ (i.e. in a polynomial ring in two variables), it is a second difficult task to recover the factorization of $P$ in the polynomial ring $\C[X_1,...,X_n]$.

\medskip

This shows the need of having an extension of Gao's Theorem for $n>2$ variables, and this will be our main result below.

\section{Hypersurface complements}

In this section $P \in \C[X_1,...,X_n]$ is a reduced polynomial
and $P={\prod_{i=1}^sP_i}$ is the factorization of $P$ into irreducible factors in $\C[X_1,...,X_n]$.
Then the associated affine hypersurface $V(P) \subset \C^n$ has $s$ irreducible components.

First we recall a basic result, relating the number $s$ of irreducible factors to the topology of the
hypersurface complement $M(P)=\C^n \setminus V(P).$

\begin{prop}\label{p1}
$$s=\dim H^1(M(P),\C).$$
\end{prop}

\proof Using Corollary 1.4 on p.103 in \cite{D1}, we get $H_1(M(P),\Z)=\Z^s$.
Then we use the usual identification $H^1(M(P),\C)=\Hom (H_1(M(P),\Z),\C)$.
\endproof

Using this result, we can give the following description of the generic 2-planes $E$.

Let $\overline{V(P)}\subset \PP^n$ be the projective
closure of the hypersurface $V(P).$ Then $E$ is said to be {\it geometrically generic with respect to $V(P)$} if its projective closure $\overline{E}$ is
 transversal to every strata of a Whitney stratification
of $\overline{V(P)}$.
Applying the Zariski Theorem of Lefschetz type, see for instance \cite{D1}, p. 25, we get the following.
\begin{cor}\label{c2}
Let $E$ be a geometrically generic affine 2-plane with respect to the affine
hypersurface $V(P)$. Then $V(P)$ and $V(P) \cap E$ have the same number of irreducible components.

\end{cor}
\proof
The Zariski Theorem of Lefschetz type implies that the two complements $M(P)$ and $E \setminus (V(P) \cap E)$
have isomorphic fundamental groups. Since we know that, for any path connected space $X$, the abelianization $ab(\pi_1(X))$ of the fundamental group coincides to the integral first homology group $H_1(X,\Z)$,
the result follows using Proposition \ref{p1}.
\endproof

For any $n$-tuple $A=(A_1,...,A_n) \in\mathbb{C}[X_1,...,X_n]^n$ of polynomials, consider the rational 1-form
$$\omega(A)=\sum_{i=1}^n\frac{A_i}{P}dX_i$$
defined  on the affine open set $M(P)$.
Such a form  $\omega(A)$ is closed by definition if
$$d\omega(A)=\sum_{i=1}^n\left[\sum_{j=1,j\neq i}^n\left(\frac{A_i}{P}\right)_{X_j}dX_j\right]\wedge dX_i=0$$
where the subscript $X_j$ means taking the partial derivative with
respect to $X_j$. In other words, the following equations should be
satisfied.
\begin{equation} \label{eq5}
\left(\frac{A_j}{P}\right)_{X_i}-\left(\frac{A_i}{P}\right)_{X_j}=0
\end{equation}
for all $i,j=1,..,n$ with $i<j$.
Consider the vector space
$F(P)$ of all solutions $A=(A_1,...,A_n)\in \mathbb{C}[X_1,...,X_n]^n$ of the equations \eqref{eq5}
with the following multi-degree bounds
$$\multideg(A_i)\leq(m_1,...,m_i-1,...,m_n)$$
for all $i=1,...n.$ Here $\multideg(P)=(m_1,...,m_i,...,m_n)$, and this obviously means
that $\deg_{X_i}P=m_i$ for all $i=1,...n$ exactly as in Definition \ref{d1}.

Any closed form $\omega(A)$ gives rise to a cohomology class $[\omega(A)] \in H^1(M(P),\C)$,
if we work with the de Rham cohomology groups of the affine smooth variety $M(P)$.

\begin{thm} \label{t2}
The linear map $T:F(P) \longrightarrow H^1(M(P),\C)$ defined by
$$T(A)=[\omega(A)]$$
is an isomorphism. In particular
$\dim F(P)=s.$
\end{thm}

\proof

To prove the surjectivity of the map $T$, we recall that a basis for the first
de Rham cohomology group $H^1(M(P),\C)$ is given by the rational 1-forms
\begin{equation} \label{eq6}
\frac{dP_j}{P_j}=\omega(B^j)
\end{equation}
for $j=1,...,s$, where $B^j=(B_1^j,...,B_n^j)$ with $B_i^j=\frac{P
\cdot (P_j)_{X_i}}{P_j}$ where the subscript $X_i$ indicates the
partial derivative with respect to $X_i$. It is clear that $B^j \in
F(P)$, which yields the surjectivity of $T$.

To prove the injectivity of $T$, assume that $T(A)=0$,i.e.
\begin{equation} \label{eq7}
\omega(A)=d\alpha
\end{equation}
for some rational function $\alpha \in \Omega^0(M).$
We can restrict to the case when $\alpha$ is a rational function in view of
Grothendieck Theorem \cite{G} saying that for an affine smooth variety
the cohomology can be computed using the regular (algebraic) de Rham complex.

It follows that $\alpha$ \ is
then a regular function of the form $\alpha=\frac{Q}{P^k}$, where
 $k\geq 0$ and $Q$ is not divisible by $P.$ Then for any index
$j\in \{1,2,...,s\}$, $\alpha$ has a pole of order $k_j \ge 0$ along the irreducible component
$V(P_j).$ Working locally in the neighborhood of a smooth point $p_j$ of $V(P_j)$, we see that $d{\alpha}$ has either a pole of order zero along
$V(P_j)$ if $k_j=0,$ or a pole of order $k_j+1$ if $k_j \geq1.$ Hence in any
case we do not get a pole of order $1$.  On the other hand, by definition, the 1-form  $\omega(A)$
has poles of order at most one along any component $V(P_j)$. The equality \eqref{eq7} is possible only if
these pole orders are all zero. This occurs only if the polynomial $P$ divides all the polynomials $A_j$
for $j=1,...,n$. But this is impossible in view of the multi-degree bounds imposed on $A_j$, unless all
$A_j$ are zero.

\endproof

\section{Finding the irreducible factors of $P$}

In this section we explain how to find the irreducible factors of $P$. Our approach is similar to that
of Gao explained in \cite{GC}, (9.2.10)-(9.2.12), but we pay more attention to a degenerate case that may occur, which explains our next definition.
\begin{defn} \label{d3}
We say that a polynomial $P \in \mathbb{C}[X_1,...,X_n]$ is
$X_1$-generic if the restriction of the projection $\pi_1:\C^n \to
\C^{n-1}$, $(x_1,x_2,...,x_n) \mapsto (x_2,...,x_n)$ to the
hypersurface $V(P)$ has finite fibers.
\end{defn}

This property, which replaces the condition  gcd $(P,P_{X_1})=1$ in Gao's approach in \cite{Gao}, can be tested by computer since we have the following obvious result.

\begin{lem} \label{l3}
Let $P=a_0X_1^m+a_1X_1^{m-1}+...+a_m$ where the coefficients $a_j$ are polynomials in $\mathbb{C}[X_2,...,X_n]$. Then $P$ is $X_1$-generic if and only if the ideal spanned by $a_0,...,a_m$
coincides to the whole ring $\mathbb{C}[X_2,...,X_n]$.
\end{lem}

\begin{ex} \label{e3} \ \\
(i) If $d$ is the total degree of $P$ and if the monomial $X_1^d$ occurs in $P$ with a non-zero coefficient,
then clearly the polynomial $P$ is $X_1$-generic. Starting with any polynomial $P$, we can arrive at this situation by making a linear coordinate change $\tilde X_1=X_1$, $\tilde X_j=X_j+c_j \cdot X_1$ for $j>1$
and suitable constants $c_j \in \C$.

(ii) Let $n=3$ and consider the polynomial $P=X^2Y^2Z^2+X$. Then $P$ is $X$ generic, but not $Y$-generic.

\end{ex}

We assume in the sequel that the polynomial $P$ is $X_1$-generic and define the following two associated vector spaces. Let $E(P)=pr_1(F(P))$, where
$$pr_1:\mathbb{C}[X_1,...,X_n]^n \to \mathbb{C}[X_1,...,X_n]$$
denotes the projection on the first factor. Let $\overline E (P)$ be the image of $E(P)$ under the canonical projection $p: \mathbb{C}[X_1,...,X_n] \to Q(P)$, where we introduce the quotient ring
$Q(P)=\mathbb{C}[X_1,...,X_n]/(P)$.

\begin{prop} \label{p5}
If the polynomial $P$ is $X_1$-generic, then the following hold.

\noindent (i) gcd $(P,P_{X_1})=1$, where the subscript $X_1$
indicates the partial derivative with respect to $X_1$.

\noindent (ii) $\dim \overline E (P)=s.$

\end{prop}

\proof To prove (i), it is enough to show that any irreducible factor $P_k$ of $P$ does not divide $P_{X_1}$. Now, with the notation from the proof of Theorem \ref{t2}, we have
\begin{equation} \label{eq8}
P_{X_1}=\sum_{j=1,s}B^j_1.
\end{equation}
In this sum, all the terms are divisible by $P_k$, except possibly
$$B^k_1=\frac{P \cdot (P_k)_{X_1}}{P_k}.$$
This term is divisible by the irreducible polynomial $P_k$ exactly when $(P_k)_{X_1}=0$ (otherwise $\deg_{X_1}P_k > \deg_{X_1}(P_k)_{X_1}$).
But $(P_k)_{X_1}=0$ implies that $P_k \in \mathbb{C}[X_2,...,X_n]$ and then, for any $b \in \C^{n-1}$ such that
$P_k(b)=0$ (which exists since $\deg P_k >0$), the line $\pi_1^{-1}(b)$ is contained in the hypersurface $V(P)$.
This contradicts the hypothesis that $P$ is $X_1$-generic, and thus proves (i).

 To prove (ii), it is enough to show that the classes of the elements $B^j_1$ for $j=1,...,s$ are linearly independent
 in $Q(P)$. Assume there is a relation
 $$\sum_{j=1,s}c_j \cdot B^j_1=C \cdot P$$
 where $c_j \in \C$ and $C \in \mathbb{C}[X_1,...,X_n]$. Checking as above the divisibility by $P_k$, it follows that
 the coefficient $c_k$ has to vanish, for all $k=1,...,s.$

 \endproof

 Exactly as in the proof above, one can show that the classes of the elements
 \begin{equation} \label{eq9}
C^j_1=B^j_1 \cdot P_{X_1}
\end{equation}
for $j=1,2,...,s$ are linearly independent in $Q(P)$. It follows that the linear subspace
 \begin{equation} \label{eq9.5}
\tilde E (P)=\{[v \cdot P_{X_1}  ]~~|~~ v \in  E (P)\}
\end{equation}
in $Q(P)$ is $s$-dimensional. Let $S: \tilde E (P) \to \overline E (P)$ be the inverse of the linear isomorphism $\overline E (P) \to \tilde E (P)$ sending
$[v]$ to $[v \cdot P_{X_1} ]$ for $j=1,2,...,s$.

 Note that in the quotient ring $Q(P)$ one has
  \begin{equation} \label{eq10}
[B^i_1]\cdot [B^j_1]=0
\end{equation}
for $i \ne j$ and
  \begin{equation} \label{eq11}
[B^i_1]\cdot [B^i_1]=[P_{X_1}]\cdot [B^i_1]
\end{equation}
for all $i=1,...,s.$ Let $v \in E(P)$ and write $[v]=\sum_{j=1,s}\lambda_j [B^j_1]$ in $Q(P)$.
Consider the linear mapping
$$\phi_v:Q(P) \to Q(P)$$
induced by the multiplication by $v$. Then the equations \eqref{eq9}, \eqref{eq10}, \eqref{eq11} imply
that $\phi _v(\overline E (P)) \subset \tilde E (P)$. It follows that  $\psi_v=S \circ \phi_v$ as a linear
endomorphism of the $s$-dimensional vector space $\overline E (P)$. We also get
$$\psi_v([B^i_1])=\lambda_i [B^i_1]$$
for all $i=1,...,s.$
  \begin{rk} \label{r2}
A key point here is that the vector space $\overline E (P)$ and the endomorphism $ \psi_v: \overline E (P) \to \overline E (P)$ can be computed without knowing the factorization of $P$.
\end{rk}
We have the following basic result.

\begin{prop} \label{p6}
If the polynomial $P$ is $X_1$-generic and all the eigenvalues of the endomorphism $ \psi_v: \overline E (P) \to \overline E (P)$ are distinct, say $\lambda_1,...,\lambda_s$, then, up-to a re-indexing of the factors, one has
$$P_i=gcd(P,v-\lambda_i P_{X_1})$$
for $i=1,...,s.$
\end{prop}

\proof

Using the above equations, we get a polynomial $C_1 \in \mathbb{C}[X_1,...,X_n]$ such that
  \begin{equation} \label{eq12}
 v-\lambda_i P_{X_1}= \sum_{j \ne i}(\lambda_j-\lambda_i) [B^j_1]+ C_1\cdot P.
\end{equation}
It follows that the irreducible polynomial $P_i$ divides $v-\lambda_i P_{X_1}$. Moreover, exactly as in the proof of Proposition \ref{p5}, we see that the irreducible polynomial $P_k$ does not divide $v-\lambda_i P_{X_1}$ for $k\ne i$.

\endproof

\end{document}